\documentclass{ifacconf}

\usepackage{graphicx}      % include this line if your document contains figures
\usepackage{natbib}        % required for bibliography
\usepackage{amsmath, amssymb, amsfonts, amscd, mathrsfs,bm}
\renewcommand{\epsilon}{\varepsilon}
\def\lip{\mathop{\mathrm{Lip}}\nolimits}

\renewcommand{\d}{\, \mathrm{d}}

\newcommand{\id}{\mathop{\bf id}}

\newcommand{\R}{\mathbb{R}}

\usepackage{algorithm}
\usepackage{algorithmic} 

\begin{document}
\begin{frontmatter}

\title{Indirect Methods in Optimal Control\\ on Banach Spaces\thanksref{footnoteinfo}} 

\thanks[footnoteinfo]{
RC and MS acknowledge financial support of the %SYSTEC – 
Research Center for Systems and Technologies 
(UID/00147) and the Associate Laboratory ARISE – %Advanced Production and Intelligent Systems 
(LA/P/0112/2020, DOI: 10.54499/LA/P/0112/2020), both funded by Fundação para a Ciência e a Tecnologia, I.P./MECI through national funds. 
A part of the computations was carried out on the OBLIVION Supercomputer (Évora University, Portugal) under FCT computational project 2023.10674.CPCA. NP is supported by a subsidy from the Ministry of Education and Science of Russia (project no.~121041300060-4).
}

\author[First]{Roman Chertovskih}
\author[Second]{Nikolay Pogodaev}
\author[First]{Maxim Staritsyn}
\author[First]{A. Pedro Aguiar}

\address[First]{Research Center for Systems and Technologies (SYSTEC),\\
Faculdade de Engenharia, Universidade do Porto, Portugal \\
(e-mails: roman@fe.up.pt, apra@fe.up.pt, starmaxmath@gmail.com)}
\address[Second]{Matrosov Institute for System Dynamics and
Control Theory (ISDCT~SB~RAS), Irkutsk, Russia (e-mail: nickpogo@gmail.com)}

\begin{abstract}                % Abstract of 50--100 words
This work focuses on indirect descent methods for optimal control problems governed by nonlinear ordinary differential equations in Banach spaces, viewed as abstract models of distributed dynamics. As a reference line, we revisit the classical schemes, rooted in Pontryagin’s maximum principle, and highlight their sensitivity to local convexity and lack of monotone convergence. We then develop an alternative method based on exact cost-increment formulas and finite-difference probes of the terminal cost. We show that our method exhibits stable monotone convergence in numerical analysis of an Amari-type neural field control problem.
\end{abstract}


\end{frontmatter}
%===============================================================================

\section{Introduction}

% Optimal control theory of evolution equations is a powerful language for PDEs and...

% Numerical algorithms are still fragmentary...

% PMP for such systems is tricky and hardly used in practice, particularly, for designing  indirect numeric methods.

% We propose a different approach based on exact increment formulas leading to so-called monotone descent methods.

% Typically, such methods are derived for state-linear models, but we learned how to do this in nonlinear frameworks. 

\subsection{Motivation}

\vspace{-0.25cm}

The language of ordinary differential equations (ODEs) in Banach spaces is natural for the analysis of many distributed systems and partial differential equations (PDEs), in particular in the context of control and optimization, see, e.g., \cite{bensoussan2007representation,Fattorini1999,Li1995} (we only mention a few monographs and refer to citations therein). The key advantage of this viewpoint is the immediate availability of ``classical'' control-theoretical tools such as Pontryagin's maximum principle (PMP). This advantage, however, is rarely exploited in numerical analysis. In particular, the framework of \emph{indirect}  methods for optimal control on Banach spaces is still rather fragmentary compared to the finite-dimensional case, see \cite{borzi2023sequential,srochko1982computational}. 
By ``indirect'' we mean numerical algorithms for optimal control that rely on principles of \emph{dynamic} extremality. In contrast to direct or semi-direct methods~--- based on a full or partial discretization of the model followed by the solution of the resulting finite-dimensional control or mathematical programming problem~--- such schemes are less sensitive to the curse of dimensionality, and their outcomes are directly interpretable in terms of the original problem. 
The goals of this work are to highlight these issues and to partially bridge the mentioned gap.

\subsection{Contribution \& Novelty}

\vspace{-0.25cm}

When talking about ODEs in infinite-dimensional spaces, one often has in mind an evolution equation with an unbounded linear operator, a generator of a $C_0$-semigroup. In this paper, we restrict ourselves to a simpler, ``pure'' ODE setting. Specifically, we focus on models representable by an ODE in a Banach space driven by a bounded (but possibly nonlinear) vector field~--- a case which still covers a wide class of relevant distributed models, yet is (in our opinion, unfairly) sidestepped in the literature. Among prominent examples are integro-differential equations with bounded kernels, delay equations, and nonlocal transport-type models with bounded interaction operators.

Most works on indirect methods for PDE-constrained optimization essentially exploit the explicit PDE structure, see \cite{Borzi2011,Troeltzsch2010}. This approach, although well established in many situations, lacks universality and typically requires substantial theoretical preparation in each concrete case. 

Passing to the ODE framework brings the PMP back to its ``elementary form'' enabling the application of existing methods 
``as they stand''. More importantly, this conversion opens a direct path to \emph{monotone methods}, i.e., algorithms based on \emph{exact} cost-increment formulas which quantify the change of the objective under a ``switch'' of the reference control to an arbitrary admissible one. In contrast to PMP-inspired schemes, such algorithms avoid model linearization and provide monotone convergence without hidden parametric search. Well known in linear problems and only recently extended to nonlinear ODEs in \emph{finite} dimensions by \cite{SChPP-2022,pogodaevExactFormulaeIncrement2024}, this approach seems to have never been established on the abstract ground of Banach space-valued dynamics. The present study is, to the best of our knowledge, the first step in this direction.

\subsection{Organization of the paper}

\vspace{-0.25cm}

The paper is organized as follows: In \S~2, we formulate the optimal control problem, specify basic regularity assumptions, and establish the existence of a solution. \S~3 recalls common PMP-based numerical algorithms. \S~4 contains the main results: we develop a monotone feedback method based on exact cost-increment formulas. This approach is then applied to an Amari-type neural field model in \S~5. All technical proofs are collected in the Appendix.

% We begin with figuring out a particular, yet sufficiently 

% formulate a general based on Pontrygin's. This approach, however, fails    

% Next, we take a different path based on \emph{exact} cost-increment formulas that quantify the change of the objective under  ``switching'' of the reference control to \emph{any} other admissible input. This approach avoids linearization of the state equation and leads to a so-called \emph{monotone} algorithm.

% In contrast to the PMP, which rests at the linear approximation of the cost functional in the class of needle-shaped (or convex) control perturbations, our method is This approach avoids linearization of the state equation and leads to monotone descent steps. % with an automatic descent property at each {\rm iter}ation. 

% Standard monotone schemes are derived for state--linear dynamics, using the duality arguments and elementary computations, see, e.g. \cite{?}. The first generalization of such algorithms and the underlying optimality principle to a general nonlinear problem was performed by \cite{pogodaevExactFormulaeIncrement2024} in the classical  setting of finite-dimensional ODEs.  In this paper, we extend the same machinery to semilinear evolution equations on Banach spaces. %, provided by the linear--quadratic structure w.r.t. the control variable. 
% % As a main result, we derive %an exact increment formula and 
% % a monotone numerical method, which is free of any intrinsic parameters,  simple to implement, and robust under the weak convergence of controls.

\section{Problem Statement}

\vspace{-0.1cm}

Let $X \doteq (X,\|\cdot\|_{X})$ be a real Banach space, $U \doteq (U,\langle\cdot,\cdot\rangle)$ a real separable Hilbert space, and let $T,R>0$ and $\alpha\geq 0$ be given parameters. On the fixed time interval $I \doteq [0,T]$ consider the optimal control problem
\begin{equation*}
({\bf P}) \quad 
\inf\Big\{\mathcal I[u] \doteq \ell(x_T^u) + \frac{\alpha}{2}\int_I\|u_t\|^{2}_{U} \d t\colon \ u \in \mathcal U\Big\}.
\end{equation*}
Here $\ell\colon X \to \mathbb R$ is a given terminal cost functional, control inputs $u\colon I \to U$, $t \mapsto u_t$, are chosen from the closed ball
\[
    \mathcal U \doteq \Big\{u \in L^\infty(I;U)\colon \|u\|_{U} \leq R \ \mbox{ for almost all } t \in I\Big\}
\]
of radius $R$ in the space of Bochner-measurable functions $I \to U$, and state trajectories $x \doteq x^u\colon I \to X$, $t \mapsto x_t$, satisfy the ODE
\begin{equation}\label{de1}
   \dot x_t = F_t(x_t,u_t) \doteq f_t(x_t) + G_t(x_t)\,u_t 
   \quad \text{for a.a. } t \in I,
\end{equation}
with a given initial condition $x_0 \in X$. The maps \[f\colon I \times X \to X, \qquad G\colon I \times X \to \mathcal L(U,X)\] are given; their regularity will be specified below.
% To keep the presentation light, we %deliberately 
% sacrifice a layer of non-essential generality. %by dropping, say, the explicit measurable dependence on time in the dynamics. 
% At the same time, the overall 
The overall structure of the problem $(\mathbf P)$ is not only dictated by technical convenience but is essential for some results of consequent sections.

\vspace{-0.15cm}

\subsection{Notation}

\vspace{-0.1cm}

Measurability of a map $\phi\colon I \to Y$ into a Banach space 
$Y \doteq (Y,\|\cdot\|_{Y})$ is understood in the strong (Bochner) sense. 
For $p \geq 1$, the Lebesgue--Bochner spaces $L^p(I;Y)$ are defined as the sets of equivalence classes of measurable functions $\phi\colon I \to Y$ such that $\big(t \mapsto \|\phi_t\|_{Y}\big) \in L^p(I)$, where two functions are identified if they coincide at almost all (a.a.) $t \in I$. By $C(I;Y)$ we denote the space of all continuous functions $\phi\colon I \to Y$; it is endowed with the uniform norm
\(
\|\phi\|_\infty \doteq \sup_{t \in I} \|\phi_t\|_{Y}.
\)
Given another Banach space $Z$, we write $\mathcal L(Y;Z)$ for the space of bounded linear operators $L\colon Y \to Z$; this space is equipped with the operator norm
\(
\|L\|_{\mathcal L(Y;Z)} \doteq \sup_{\|\mathrm y\|_{Y} \leq 1} \|L \mathrm y\|_{Z}.
\)
$L'\colon Z' \to Y'$ is the adjoint of $L$. The class $C^1(Y;Z)$ is formed by maps $F \colon Y \to Z$ which are continuously differentiable in the Fr\'echet sense, i.e. such that
\(
D F\colon Y \to \mathcal L(Y;Z)\) %,
is continuous; we omit mentioning $Z$ when it coincides with $\R$. A Hilbert space $H$ is identified with its dual $H'$ via  the Riesz isometry.

\vspace{-0.1cm}
 
\subsection{Standing Assumptions \& Preliminaries}\label{ssec:assumptio}

\vspace{-0.1cm}

The first group of regularity assumptions is standard: 
\begin{description}

\item[(${\bf H}_{1}$)] For each \(\mathrm{x}\in X\), the maps \(t\mapsto f_t(\mathrm{x})\) and \(t\mapsto G_t(\mathrm{x})\) are measurable, and there exist constants \(M_f, M_{G}\geq 0\) such that, for any ${\rm x, y} \in X$ and a.a. $t \in I$, the following conditions hold:
\[
\begin{aligned}
&\|f_t(\mathrm{x})-f_t(\mathrm{y})\|_X\le M_{f}\,\|\mathrm{x}-\mathrm{y}\|_X, \ \ 
\|f_t(0)\|_X\le M_f,\\
% \end{aligned}
% \]
% and
% \[
% \begin{aligned}
&\|G_t(\mathrm{x})-G_t(\mathrm{y})\|_{\mathcal L(U;X)}\le M_{G}\,\|\mathrm{x}-\mathrm{y}\|_X,\\ &\|G_t(0)\|_{\mathcal L(U;X)}\le M_G.
\end{aligned}
\]
\end{description}
These conditions are not sharp; we %deliberately 
sacrifice a layer of non-essential generality to keep the presentation light.

Fix $u\in \mathcal U$, $\mathrm{x}\in X$ and \(s \in [0, T)\). By a solution to the ODE \eqref{de1} with initial data $(s, \mathrm x)$ we mean an absolutely continuous function $t \mapsto \Phi^u_{s,t}(\mathrm x)$ satisfying the following Bochner integral equation
\begin{equation}
\begin{aligned}\label{Phi_integral_repr}
\Phi^u_{s,t}(\mathrm x) & = \mathrm{x}+ \int_s^t  F_\tau(\Phi^u_{s,\tau}(\mathrm x), u_\tau)\d \tau{,} \qquad t \in [s, T].
\end{aligned}
\end{equation}
For any triple $(s, \mathrm x, u) \in [0,T)\times X \times \mathcal U$, hypotheses $({\bf H}_1)$ guarantee the existence of a unique solution $\Phi^u_{s,\cdot}(\mathrm x) \in C([s,T];X)$ on the whole interval $[s, T]$; moreover, this solution is globally bounded %and Lipschitz in $\mathrm x$ on bounded subsets of $X$ 
\cite[Th.~2.2.2]{kolokoltsov2019}.

Denoting by $\id_{X}$ the identity map on $X$, it is easy to verify the following relations for any \(0\leq t_0 \leq t_1 \leq t_2 \leq T\):
\begin{align}
\label{fami}
\Phi^u_{t_1,t_2}\circ\Phi^u_{t_0,t_1}=\Phi^u_{t_0,t_2},\qquad \Phi^u_{t_0,t_0}={\id}_{X}.
\end{align}
A map $(s, t, \mathrm x) \mapsto \Phi_{s,t}^u(\mathrm x)$ satisfying such algebraic axioms is called an evolution map or a (semi-) \emph{flow} on $X$. Thanks to the time-reversibility of the ODE \eqref{Phi_integral_repr}, all $\Phi_{s,t}$ are homeomorphisms $X \to X$. Moreover, if the functions $f_t$ and \(G_t\) are $C^1$ for a.a. $t$, then $\Phi_{s,t}$ are $C^1$-diffeomorphisms, see \cite{cartan1983differential}.

\vspace{-0.15cm}

\subsection{Existence of a Minimizer}\label{parag:exist}

\vspace{-0.1cm}

In this paragraph, we specify additional regularity assumptions to ensure well-posedness of the problem $(\mathbf P)$. The argument is standard: we need to guarantee the compactness of $\mathcal U$ in a suitable weak topology and the corresponding lower semicontinuity of $\mathcal I$.

The first part is straightforward. As a Hilbert space, $U$ has the Radon--Nikodym property. Hence, $L^{\infty}(I;U)$ is identified isometrically with $(L^1(I;U))'$ via the pairing
\[
\langle\!\langle u,v\rangle\!\rangle \doteq \int_I \langle u_t,v_t\rangle \d t, \ \  u \in L^\infty(I; U), \, v \in L^1(I;U),
\]
see, e.g. \cite[\S~II]{diestel1977vector}. Endowing the closed ball $\mathcal U$ of $L^\infty$ with the weak* topology $\sigma(L^\infty, L^1)$, induced by this duality, makes the resulting topological space $\mathcal U_{w*}$ compact thanks to the classical Banach-Alaoglu theorem.

% In this paragraph, we prove that variational problem $(\bf P)$ has a minimizer in the class $\mathcal U$. Our strategy is standard and based on the application of the Weierstrass theorem, provided by the compactness of the set $\mathcal U$ in an appropriate topology, and the lower semicontinuity of the objective function $\mathcal I[\cdot]$. 

The second ingredient is given by the following lemma, whose proof is outlined in Appendix A.
\begin{lem}
In addition to $(\mathbf H_1)$, suppose the following
\begin{description}
\item[(${\bf H}_{2}$)] \(\ell\) is lower semicontinuous.
\item[(${\bf H}_{3}$)] The operator $G$ has a finite--``nuclear'' structure:
\begin{align}
\label{G}
G_t(\mathrm x) \, \mathrm u = \sum_{j=1}^m \langle \mathrm u, g^j_t(\mathrm x)\rangle \, h^j_t(\mathrm x)
\end{align}
with $g^j\colon I \times X \to U$ and $h^j\!: I \times X \to X$ fitting the related assumptions in $(\mathbf H_1)$. 
\end{description}
Then, 1) the input-output operator $u \mapsto x^u$ is continuous as a map from $\mathcal U_{w*}$ to $(C(I;X), \|\cdot\|_\infty)$, and 2) the mapping $u \mapsto \|u\|^2_{L^2(I;U)}$ is lower semicontinuous on  $\mathcal U_{w*}$. In particular, $\mathcal I$ is lower semicontinuous as $\mathcal U_{w*} \to \R$. 
\end{lem}

Applying the Weierstrass theorem, we finally obtain:
\begin{thm}
    Under hypotheses $({\bf H}_{1})$--$({\bf H}_{3})$, problem $(\mathbf P)$ admits a minimizer in the class $\mathcal U$.
\end{thm}
Hypothesis $({\bf H}_3)$ means that the system is actuated through a finite number of effective directions with scalar ``channel'' outputs. %, playing the role of aggregated control signals. 
This setting is natural for distributed-parameter systems driven by a finite family of actuators or sensors. In this case, the operators $G_t(\mathrm x)$, $\mathrm x \in X$, have finite rank, and theirs adjoints $G_t(\mathrm x)'\colon X' \to U$ are \[G_t(\mathrm x)'p = \sum \langle p, h^j_t(\mathrm x)\rangle\, g^j_t(\mathrm x),\ \ \ p \in X'.\]
%Assumptions $({\bf H}_{1})$--$({\bf H}_{3})$ remain in force for the rest of the paper. 

\vspace{-0.1cm}

\section{PMP-based Descent Methods}\label{chap:Kry}

\vspace{-0.168cm}

%\subsection{Pontryagin's minimum principle}

The PMP is a standard route to local optimality in dynamic contexts. Although it is technically delicate for general constrained evolution equations \cite[Exam.~1.4]{Li1995}, in our  case it does not principally differ from its classical counterpart of \cite{Pontryagin1962}.

In addition to $(\bf{H}_{1})$--$(\bf{H}_{3})$, assume:
    \begin{description}
        \item[$({\bf H}_{4})$] $\ell \in C^1(X)$, $f, h^j \in C^1(X;X)$,  and $g^j \in C^1(X;U)$.
    \end{description}

Given $\bar u \in \mathcal U$, denote by $\bar x$ the corresponding trajectory and by $\bar p$ the solution of the \emph{adjoint} (linear backward) ODE
\begin{equation}\label{psi}
\dot \psi_t = - D F_{t} (\bar x_t, \bar u_t)' \, \psi_t, \quad \psi_T = D\ell(\bar x_T),
\end{equation}
and introduce the Hamilton-Pontryagin functional: 
\[
H_t(\mathrm x, \mathrm p, \mathrm u) = \frac{\alpha}{2}\|\mathrm u\|^2_U + \langle  \mathrm p, f_t(\mathrm x)\rangle_{(X', X)} + \langle\mathrm u, G_t(\mathrm x)' \mathrm p \rangle.
\]
The optimality of $\bar u$ in $(\bf P)$ implies \cite[Th. 1.6]{Li1995} that
\begin{align}
H(\bar x_t, \bar \psi_t, \bar u_t) = \min_{{\rm u} \in U} H(\bar x_t, \bar \psi_t, \mathrm u) \quad  \text{ for a.a. } \ \ t \in I.\label{PMP}
\end{align}
A control process satisfying the PMP is termed \emph{extremal}. Unfortunately, the extremality is \emph{in}sufficient for the global optimum even in the finite-dimensional case. 

For $\alpha >0$ and $R$ sufficiently large, the minimizer in \eqref{PMP} is uniquely defined as $u_t^* = -\alpha^{-1} G(\bar x_t)' \bar \psi_t$. In that case, we can approach an extremal by adopting the simplest method of \cite{Krylov1963}.  Note however, that the convergence of this method is ensured under fairly strong local convexification (\(\alpha>\!\!> 1\)), resulting in high amplitudes of \(u^\ast\), which is often unphysical and affects numerical stability of the state equation. Moreover, the descent fails, in general, to be monotone. The monotonicity can be enforced by localizing the effect of \(u^\ast\) through 
a tradeoff with a baseline control \(\bar u\), as in \cite{srochko1982computational}, or by sequential quadratic Hamiltonian approximations, as in 
\cite{borzi2023sequential}. Both approaches assume a parametric line search at each descent step, typically implemented by backtracking and leading either to an essential growth in the number of iterations or to multiple recomputations of primal and adjoint states. These ideas are captured at a conceptual level by Algorithm~\ref{alg0}.

\begin{algorithm}[h!]
  \caption{(PMP)}
  \label{alg0}
  \begin{algorithmic}[]
    \REQUIRE initial guess \(u^0 \in \mathcal U\)
    \ENSURE sequence \((u^{{\rm iter}}) \subset \mathcal U\) with \(\mathcal I[u^{{\rm iter}+1}] \leq \mathcal I[u^{{\rm iter}}]\)
    \FOR{\({\rm iter} = 0,1,2,\dots\)}
      \STATE Set \(\bar u \gets u^{{\rm iter}}\)
      \STATE Solve \eqref{de1} with \(u = \bar u\) forward in time to obtain \(\bar x\)
      \STATE Solve \eqref{psi} backward in time to obtain \(\bar p\)
      \STATE Update $u^{{\rm iter}+1}_t \gets (1-\eta) \bar u_t -\eta \, \alpha^{-1} \, G(\bar x_t)' \, \bar \psi_t$\\ with $\eta \in (0,1)$ chosen dynamically (backtracking)
    \ENDFOR
  \end{algorithmic}
\end{algorithm}

Below, we propose a different scheme that achieves monotonicity through a feedback mechanism, at the price of an a priori chosen number of forward solutions to \eqref{de1}.

\section{Monotone Methods}

\vspace{-0.3cm}

The approach, we are going to develop, has a fragrance of dynamic programming but remains local in nature. The key is that the nonlinear Hamilton--Jacobi equation does not appear. Instead, we work with a linear transport equation whose solution has an explicit form. 
%The tradeoff is that the nonlinear Hamilton-Jacobi equation will not come on scene. Instead, we deal with the linear transport equation whose solution admits ex explicit representation.

Let $\bar u\in\mathcal U$ be a baseline control, $\bar\Phi\doteq\Phi^{\bar u}$ the corresponding flow, and $\bar x_t\doteq\bar\Phi_{0,t}(\mathrm x_0)$. Taking another control $u\in\mathcal U$ and a moment $s\in(0,T)$, we define an element $u\triangleright_s\bar u \in \mathcal U$ as
\begin{equation}
t\mapsto(u\triangleright_s\bar u)_t\doteq
\begin{cases}
u_t,&t\in[0,s),\\
\bar u_t,&t\in[s,T].
\end{cases}
\label{control-var-gen}
\end{equation}
%This is just a new admissible control, obtained by switching from $u$ to $\bar u$ at the time moment $s$. 
%It can be viewed as a kind of needle-shaped perturbation with a ``large support'', and customary appears in the context of dynamic programming, e.g., in the analysis of the problem value function. 
By construction, we have:
\[
x^{u\triangleright_s\bar u}_t
=\bigl(\bar\Phi_{s,t}\circ\Phi_{0,s}\bigr)(\mathrm x_0)
=\bar\Phi_{s,t}(x_s)
\]
for \(0 \leq s \leq t \leq T\), where $\Phi\doteq\Phi^{u}$ is the flow generated by $u$, and $x_t\doteq x^u_t \doteq \Phi_{0,t}(\mathrm x_0)$.  

Denote \(\bar p_t \doteq \ell\circ \bar \Phi_{t, T}\); the value \(\bar p_t(\mathrm x)\) is the terminal cost of the baseline control for the modified problem \(({\bf P}_{t,\mathrm x})\), in which the initial data \((0,\mathrm x_0)\) are replaced by \((t,\mathrm x)\). For any \(\mathrm x\), the composition property~\eqref{fami} yields:
% Denote $\bar p_t \doteq \ell\circ \bar \Phi_{t, T}$; the value $\bar p_t(\mathrm x)$ is the terminal cost of the baseline control in the modification $({\bf P}_{t, \rm x})$ of problem $(\bf P)$, obtained by taking $(t, {\rm x})$ as initial data. For any $\rm x $, the composition property \eqref{fami} gives:
\begin{gather}\label{ppp}
\bar p_t\bigl(\bar \Phi_{0, t}{(\rm x)\bigl) \, \equiv \,} \ell(\bar \Phi_{0, T}(\rm x)),
\end{gather}
saying us that the quantity on the left is independent of $t$ and, in particular, $\ell(\bar x_T) = \bar p_T(\bar x_T) = \bar p_0({\rm x_0})$. This enables representing the increment in the terminal cost rate on the pair $(\bar u, u)$ as
\begin{align}
&\ell(x_T) - \ell(\bar x_T) = \bar p_T(x_T) - \bar p_T(\bar x_T)\nonumber\\
& \ \ \ \ \ \   = \bar p_T(x_T) - \bar p_0({\rm x_0}) = \int_I \frac{\d}{\d t} \bar p_t(x_t) \d t.\label{ppl}
\end{align}
The subsequent analysis relies on computing the derivative under the sign of the integral. This is straightforward when \(\bar u\), \(f_\cdot\), and \(G_\cdot\) are continuous (so that \(t \mapsto \bar p_t\bigl(\bar \Phi_{0,t}(\mathrm x)\bigr)\) belongs to \(C^1(I;\mathbb R)\)). In such a case, differentiating the identity \eqref{ppp} and applying the standard chain rule gives:
\[
\Bigl(\partial_t \bar p_t(\mathrm y) + D\bar p_t(\mathrm y)[F_t(\mathrm y, \bar u_t)]\Bigr)\Big|_{{\rm y = \bar \Phi}_{0, t}(\mathrm x)} = 0. 
\]
Leveraging the fact that \(\bar \Phi_{0,t}\), $t \in I$, are homeomorphisms, we see that \(\bar p\) satisfies the linear PDE in \(I \times X\):
\begin{gather}
\partial_t \bar p_t + D\bar p_t[F_t(\cdot, \bar u_t)] = 0. \label{eq_p}
\end{gather}
In particular,
\(
\partial_t \bar p_t(x_t) = - D\bar p_t(x_t)[F_t(x_t, \bar u_t)].
\)
At the same time, the usual chain rule gives
\[
    \frac{\d}{\d t} \bar p_t(x_t)
    = \partial_t \bar p_t(x_t) + D\bar p_t(x_t)[F_t(x_t, u_t)].
\]
Combining these identities and using the structure of $F$, we arrive at the following key representation:
\begin{equation}
    \frac{\d}{\d t} \bar p_t(x_t)
    = \left\langle u_t - \bar u_t, G_t(x_t)' D\bar p_t(x_t) \right\rangle.
    \label{eif}
\end{equation}
A rigorous derivation of this equality in the general case requires extra regularity assumptions, which are collected in Lemma~3; its proof is given in Appendix~B.

\begin{lem}
    Along with $(\bf{H}_{1})$--$(\bf{H}_{4})$, impose the hypothesis (for brevity, in terms of the aggregated operator $G$):
    \begin{description}

        \item[$({\bf H}_{5})$] For any compact subset $K \subset X$, there exist constants $M_K(D\ell), M_K(Df), M_K(DG)\geq 0$ such that \(\displaystyle\big\|D\ell(\mathrm x)\bigl\|_{\mathcal L(X;\R)} \leq M_K(D\ell)\), \(  
        \bigl\|Df_t(\mathrm x)\big\|_{\mathcal L(X;X)} \leq M_K(Df),\)
        and 
        $
        \big\|DG_t(\mathrm x)\bigr\|_{\mathcal L(X;\mathcal L(U;X))} \leq M_K(DG)
        $ for all $\mathrm x\in K$, and a.e. $t \in I$.
    \end{description}
    Then, \eqref{eif} holds a.e. on $I$. 
\end{lem} 

By substituting \eqref{eif} inside \eqref{ppl}, we obtain an exact increment formula for the cost functional $\mathcal I$:
\[
\begin{gathered}
\mathcal I[u] - \mathcal I[\bar u]
= \int_I \Bigl(\bar H_t(x_t, u_t) - \bar H_t(x_t, \bar u_t)\Big)\d t,\\
% \]
% where %\(\bar H_t\) is a contraction of $H$ to  $\mathrm p = D\bar p_t(\mathrm x)$: 
% \[
\bar H_t(\mathrm x, \mathrm u) \doteq H_t (\mathrm x,  D\bar p_t(\mathrm x), \mathrm u)
= \frac{\alpha}{2}\|\mathrm u\|_U^2
  + \langle \mathrm u, G(\mathrm x)' D\bar p_t(\mathrm x)\rangle.
\end{gathered}
\]
For any \((t,\mathrm x)\), the minimizer of the function \(\bar H_t(\mathrm x, \cdot)\) over  the ball \(B_R \doteq \{{\rm u} \in U\colon \|{\mathrm u}\|_U \leq R\}\) is given by
\[
    \hat u_t(\mathrm x) = \begin{cases}\Pi_{B_R}\left(w_t({\rm x})\right), & \alpha >0\\
    -R\frac{G(\mathrm x)' D\bar p_t(\mathrm x)}{\|G(\mathrm x)' D\bar p_t(\mathrm x)\|_U}, & \alpha =0,\end{cases}
\] 
where \(w_t({\rm x}) = - \frac{1}{\alpha}\,G(\mathrm x)' D\bar p_t(\mathrm x)\) and $\Pi_{B_R}$ denotes the projection operator onto the set $B_R$.
%Evidently, it is unique as soon as $\alpha>0$. 
% \begin{gather*}
% \hat u_t(\mathrm x)
% = \lambda(\mathrm x) \, w_t(\mathrm x), \quad w_t(\mathrm x) \doteq  - \frac{1}{\alpha}\,G(\mathrm x)' D\bar p_t(\mathrm x),\\ \lambda(\mathrm x) \doteq \min\Big\{1, \frac{R}{\|w_\cdot(\mathrm x)\|_{L^2(I;U)}}\Big\}.
% \end{gather*}
By taking \(u = \hat u\circ x\) and \(x=x^u\), we formally obtain:
\[
\bar H_t(x_t, u_t)
\le \bar H_t(x_t, \bar u_t)
\quad\text{for a.a. } t \in I,
\]
and $\mathcal I[u] \le \mathcal I[\bar u].$ Note, however, that in this case the trajectory \(x\) is generated by the control \(u\) itself, so the law \(u_t = \hat u_t(x_t)\) defines a feedback loop. %This fact makes a substantial difference with the PMP setting.
The existence of a backfed solution \(x^u\) can be shown, e.g., by Kakutani's fixed-point argument; when \(\alpha > 0\) and $R$ are such that \(\hat u = w\) is single-valued and continuous, the fact follows directly from the Cauchy theory for the ODE \(\dot x_t = F_t(x_t, w_t(x_t))\).

The direct computation of \(D\bar p_t\) assumes solving the PDE~\eqref{eq_p} on the underlying Banach space \(X\), which is rarely feasible in practice. Instead, we can exploit the explicit representation $\bar p_t \doteq \ell\circ \bar \Phi_{t, T}$ and hypothesis $(\mathbf H_3)$. 

To fix ideas, let $m=1$, $\alpha>0$, and $R>\!\!> 1$. Then,
\[
G(\mathrm x)' D\bar p_t(\mathrm x) =  D \bar p_t(\mathrm x)[h(\mathrm x)] \, g(\mathrm x),
\]
and the term $D\bar p_t(\mathrm x)[h(\mathrm x)]$ can be approximate by 
\[
\xi^\epsilon(t,\mathrm x)
\doteq
\frac{
  \ell\bigl(\bar\Phi_{t,T}(\mathrm x + \epsilon \, h(\mathrm x))\bigr)
  - \ell\bigl(\bar\Phi_{t,T}(\mathrm x)\bigr)
}{\epsilon}
\]
with a sufficiently small  $\epsilon>0$.
% In an orthonormal basis of \(U\) adapted to the channels \(g^j\) we may interpret
% \((\xi^h_1(t,\mathrm x),\dots,\xi^h_m(t,\mathrm x))\) as an approximation of the vector \(G(\mathrm x)' D\bar p_t(\mathrm x)\).
We now discretise the time interval \(I\) and construct an approximate feedback in a sample-and-hold fashion: Let
\(
0 = t_0 < t_1 < \dots < t_N = T
\)
be a (uniform) partition of \(I\) into \(N \ge 1\) subintervals and \(\mathrm x^k \doteq x_{t_k}\) denote the states of the controlled system at the grid points. At each step \(k\), we freeze \(\mathrm x^k\) and approximate the feedback law $\hat u$ on each subinterval \([t_k,t_{k+1})\) by 
\[
\widetilde u_t \doteq  - \alpha^{-1}\xi^\epsilon(t,\mathrm x^k) \, g(\mathrm x^k).
\]
The ODE \eqref{de1} is then integrated on 
\([t_k,t_{k+1}]\) with the initial condition 
\(x_{t_k} = \mathrm x^k\) and control 
\(\widetilde u\), so that the next sample state is defined as \(\mathrm x^{k+1} = x_{t_{k+1}}\). Iterating over \(k = 0,\dots,N-1\) produces the process \((x, u)\) on the whole interval \(I\). This synthesis mechanism, which mimics the approach of \cite{krasovskii2011game}, is summarized in Algorithm~\ref{alg:sample_hold}. For \(\varepsilon \downarrow 0\) and \(N \uparrow \infty\) (the parameters can be unified by choosing $\epsilon N = 1$), this algorithm produces a monotonically nonincreasing (and so, \emph{convergent}) sequence \(\bigl(\mathcal I[u^{\rm iter}]\bigr)\). Moreover, it can be shown that any partial limit of \(\bigl(u^{\rm iter}\bigr)\) is a ``feedback extremal'' in the sense of \cite{pogodaevExactFormulaeIncrement2024}.

\begin{algorithm}[h]
  \caption{(Monotone Descent)}
  \label{alg:sample_hold}
  \begin{algorithmic}[]
    \REQUIRE initial guess \(u^0 \in \mathcal U\), \(N \geq 1\) (\# sample points), \(\epsilon\,{>}\,0\) ( probe radius)
    \ENSURE sequence \((u^{{\rm iter}})\subset \mathcal U\) with \(\mathcal I[u^{{\rm iter}+1}] < \mathcal I[u^{{\rm iter}}]\)
    \FOR{\({\rm iter} = 0,1,2,\dots\)}
    \STATE Set \(\bar u \gets u^{{\rm iter}}\),  \(\mathrm x^0 \gets \mathrm x_0\)
    \FOR{\(k = 0,1,\dots,N-1\)}
 % \FOR{k = 0,1,\dots,N-1}
  \STATE Compute $\xi^\varepsilon(\cdot,\mathrm x^k)$ using the baseline control $\bar u$
  \STATE Define $\widetilde u_t := -\alpha^{-1}\,\xi^\varepsilon(t,\mathrm x^k)\,g(\mathrm x^k)$ for $t\in[t_k,t_{k+1})$
  \STATE Solve \eqref{de1} on $[t_k,t_{k+1}]$ with $x_{t_k} = \mathrm x^k$ and $u = \widetilde u$
    \ENDFOR
    \STATE Update $u^{{\rm iter}+1} \gets u$
    \ENDFOR
  \end{algorithmic}
\end{algorithm}

%~\ref{alg:sample_hold} (Apendix B). %  that approximates the continuous extremal feedback \(u_t = w_t(x_t)\).

\vspace{-0.1cm}

\section{Case Study}

\vspace{-0.2cm}

We now illustrate the abstract framework with a simple yet substantive model from mathematical neuroscience; such neural field equations date back to the seminal work of \cite{amariDynamicsPatternFormation1977}, but have recently gained renewed interest in connection with modern studies in machine learning:
\begin{gather}
  \partial_t N(t,\theta)
  = -\gamma\,N(t,\theta) + (W*\sigma(N(t,\cdot)))(\theta) + S(t,\theta)\notag,\\
  N(0,\theta)=N_0(\theta), \quad \theta\in[0,2\pi).\label{ex}
\end{gather}
Here, \(N(\cdot,\theta)\) represents the mean activity of a population of neurons indexed by \(\theta\), viewed as a feature variable on the unit circle \(\mathbb S^1\). The dynamics describe the balance between passive decay with rate \(\gamma>0\), recurrent synaptic input mediated by the connectivity kernel \(W\), nonlinear firing-rate response \(\sigma\), and an external control field \(S\). In the lateral-inhibition setting, the kernel encodes short-range excitation and longer-range inhibition, leading to the formation of localized ``bumps'' of activity. 

The synaptic kernel is chosen in the von~Mises form
\[
  W(\Delta)
  = \frac{e^{\kappa\cos\Delta}}{2\pi\,I_0(\kappa)},\ \, \Delta=\theta-\theta', \ \   I_0(\kappa)
  = \frac{1}{2\pi}\int_{0}^{2\pi} e^{\kappa\cos\theta}\,\mathrm d\theta,
\]
where \(\kappa>0\) is a concentration parameter, \(I_0\) the modified Bessel function of the first kind (order~0), and ``\(*\)'' the convolution on \(\mathbb S^1\). The nonlinearity
\[
  \sigma(q)=\frac{1}{1+e^{-\beta(q-\vartheta)}}
\]
is a logistic activation with slope \(\beta>0\) and threshold \(\vartheta\in\mathbb R\). The field $S$ is decomposed in a truncated orthonormal Fourier basis in \(\theta\) with coefficients \(u \doteq (u_0,u^c_k,u^s_k)_{k=1}^K\):
\[
  S(t,\theta)=u_0(t)\,\phi_0(\theta)
  + \sum_{k=1}^{K}\Bigl(u^{c}_k(t)\,\phi^{c}_k(\theta)
  +u^{s}_k(t)\,\phi^{s}_k(\theta)\Bigr),
\]
where
\(
  \phi_0(\theta)=\frac{1}{\sqrt{2\pi}}\), \(\phi^{c}_k(\theta)=\frac{\cos(k\theta)}{\sqrt{\pi}}\), \(\phi^{s}_k(\theta)=\frac{\sin(k\theta)}{\sqrt{\pi}},\)
 \(\langle \phi_i,\phi_j\rangle_{L^2(0,2\pi)}=\delta_{ij}\); the components of \(u\) are amplitudes of spatially structured stimuli serving as actual controls.

We now embed~\eqref{ex} into the abstract control system~\eqref{de1}. Let
$X \doteq C(\mathbb S^1)$
be the Banach space of continuous \(2\pi\)-periodic functions \(x\colon[0,2\pi]\to\mathbb R\) with the supremum norm \(\|\cdot\|_X=\|\cdot\|_\infty\), and
\(
  U \doteq \mathbb R^{2K+1}
\)
with the standard Euclidean inner product, so \(u_t\in U\) collects the Fourier coefficients of the control at time \(t\). Defining
\[
  K(\mathrm y)(\theta)\doteq (W * \mathrm y)(\theta)
  = \int_0^{2\pi} W(\theta-\theta')\,\mathrm y(\theta')\,\mathrm d\theta',
\]
the PDE~\eqref{ex} can be written in the form~\eqref{de1} with
\begin{gather*}
  f(\mathrm x)(\theta)
  \doteq -\gamma\,\mathrm x(\theta)
        + \bigl(K(\sigma(\mathrm x))\bigr)(\theta),\\
  (G \mathrm u)(\theta)
  = \mathrm u_0\,\phi_0(\theta)
  + \sum_{k=1}^{K}\Bigl(\mathrm u^{c}_k\,\phi^{c}_k(\theta)
  + \mathrm u^{s}_k\,\phi^{s}_k(\theta)\Bigr).
\end{gather*}
Thus, the control dependence is affine and realizes Hypothesis $({\bf H}_3)$ with a fixed finite family of spatial directions \(\{\phi_0,\phi^c_k,\phi^s_k\}\) and constant channel gains. The remaining assumptions $(\mathbf H_1)$--$(\mathbf H_2)$ and $(\mathbf H_4)$--$(\mathbf H_5)$ are easy to check.

% One checks that all assumptions $(\mathbf H_1)$--$(\mathbf H_4)$: the convolution operator \(K\) is bounded on $C(\mathbb S^1)$, the activation $\sigma$ is smooth and Lipschitz on bounded subsets of $\mathbb R$, and both $f$ and $G$ are continuously differentiable with locally bounded derivatives on $X$. Therefore all abstract results of Sections~\ref{parag:exist} and~\ref{chap:Kry} are applicable to this setup.

Consider the problem of steering the terminal state $x_T \doteq N(T,\cdot)$ towards a prescribed bump-like target profile
\[
  N_{\mathrm des}(\theta)
  = \frac{A_d\,e^{\kappa_d\cos(\theta-\theta^\ast)}}{2\pi\,I_0(\kappa_d)},
\]
which models a localized activity pattern centered at \(\theta^\ast\) with amplitude \(A_d\) and concentration \(\kappa_d\). The performance index penalizes the terminal tracking error together with the $L^2$-energy of the control:
\[
  \mathcal I[u]
  = \frac12\int_{0}^{2\pi}\!\bigl(N(T,\theta)-N_{\mathrm des}(\theta)\bigr)^2\,\mathrm d\theta
  + \frac{\alpha}{2}\|u\|_{L^2(I;U)}^2.
\]
The algorithms, developed in the previous section, are implementable in terms of the Fourier coefficients \(u_0,u_k^c,u_k^s\). 
Both algorithms used the pseudospectral method, where the derivatives and convolutions are computed in the spectral (Fourier) space, and multiplications in the original (state) space; we used fast Fourier transforms in $\theta$ to switch between these spaces. Time integration was performed by the standard 4th order Runge-Kutta method with a constant time step.

%\section{Numerical Experiments}

\begin{figure}
    \centering
    \includegraphics[width=0.62\linewidth]{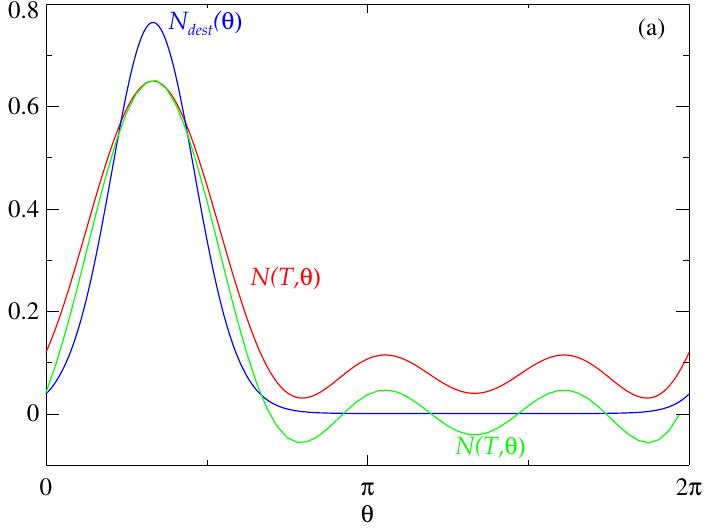}
    \includegraphics[width=0.63\linewidth]{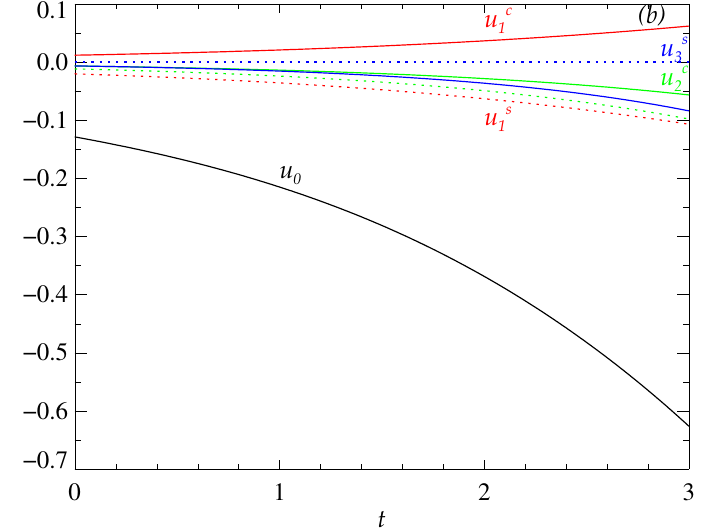} \includegraphics[width=0.63\linewidth]{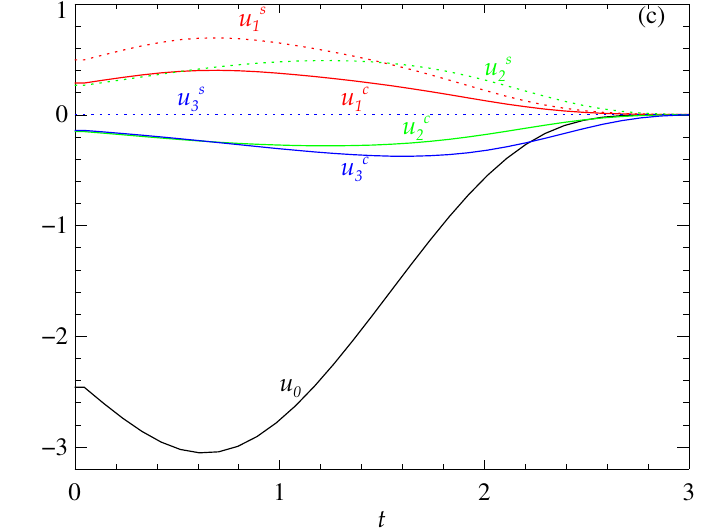}
    \caption{\small Panel (a): terminal distribution  $N(T,\cdot)$, optimized by Alg.~1 (red) and Alg.~2 (green) vs. target profile $N_{\mathrm des}$ (blue). The controls generated by Alg.~1 (b) and Alg.~2 (c). % as a function of time, $u_k^c$ are plotted by solid and $u_k^c$ dotted lines of the same color for same value of $k$.
    }
    \label{fig:placeholder}
\end{figure}

For simulation, we have chosen the following parameters: \(T=3; \ N_0(\theta)\equiv 0; \ \gamma=1,\ \beta = 2.0, \
  \kappa=4.0,\
  A_d=0.8,\ \kappa_d=6.0,\ \theta^\ast=\tfrac{\pi}{3}, \ K =3\) and 
$\alpha = 0.1$ (so $\frac{\alpha}{2}\|u\|_{L^2(I;U)}^2$ serves as a regularizing term rather than physical payoff). The outputs of Algorithms~\ref{alg0} and~\ref{alg:sample_hold}, both initialized by %equilibrium profile 
\(u^0 \equiv 0\), are shown in Fig.~\ref{fig:placeholder}. For Alg.~\ref{alg0}, we performed \(40\) iterations with an adjustable tradeoff \(\eta\). Alg.~\ref{alg:sample_hold} attained a comparable result in a \emph{single} iteration with \(N = 32\) control switchings (smoothed a posteriori). Overall, Alg{.}~\ref{alg:sample_hold} is substantially faster and yields a more accurate matching of the target mean. 
Notably, the optimized controls look different: pure exponents in the first case, and harmonic functions in the second. This could be the evidence of different extremals and deserves further investigation.  

\section{Conclusions}

Although the proposed approach performs well in experiments, its main benefit, in our view, is due to theoretical transparency and direct links to the theory of feedback extremality, which remained largely outside the scope of this short paper. A systematic comparison with existing indirect schemes and modern direct or hybrid ``discretize--optimize'' methods would require a dedicated benchmarking study. Likewise, a full, in particular quantitative, convergence analysis of the algorithm is left for future works.

\appendix

\section{Proof of Lemma 1}

(1) The continuity of the input-output operator $u \mapsto x^u$, $\mathcal U_{w*}\to (C(I;X), \|\cdot\|_\infty)$ follows from standard arguments based on Gr\"onwall's inequality. These are simple but technical and will be
omitted for brevity. 

(2) Let us ensure the lower semicontinuity of the integral term in $\mathcal I$. Since $I$ has finite measure, any $u \in \mathcal U$ belongs to $L^2(I;U)$, so the integral term can be seen as the map $u \mapsto \|u\|^2_{L^2(I;U)}$ from the Hilbert space $L^2(I;U)$. This map is weakly lower semicontinuous (as well as any strongly continuous convex functional on a Hilbert space). Finally, the inclusion $L^2 \subset L^1$ implies that any sequence converging weakly* in $L^\infty$ converges weakly in $L^2$ to the same limit. This brings the desired result.

\vspace{-0.13cm}

\section{Proof of Lemma~3}

\vspace{-0.14cm}
We use $G$ as a shorthand of \eqref{G} and drop the explicit dependence of $f$ and $G$ on $t$ for simplicity.

First, choose any representatives of the classes $\bar u$ and $u$ to be, with a slight abuse of notation, denoted by the same symbols. Recall that the map $t \mapsto %\bar{p}_t() \doteq 
\ell(\bar \Phi_{t, T}(x_{t}))$ is Lipschitz, and therefore, it is a.e. differentiable on $I$ by Rademacher's theorem. Let $J$ denote the set of differentiability points, where both $\bar u_t$ and $u_t$ are defined and belong to $U$. By construction, $J$ has Lebesgue measure of the whole $I$. 

Fixing an arbitrary $t \in J \cap [0, T)$, we need to show that
\[
\lim_{h\downarrow 0}\frac{\bar p_{t+h}(x_{t+h}) - \bar p_t(x_{t}) - h D \bar p_t(x_t)\big[G(x_t)[u_t - \bar u_t]\big]}{h} = 0. 
\]
To this end, we perform the following steps.

1) \emph{Decomposition of the limit and uniform estimates.} For $h \in (t, T]$, we represent:
\[
\begin{aligned}
&\bar p_{t+h}(x_{t+h}) - \bar p_t(x_{t}) =\\ 
&\qquad \bar p_{t+h}(x_{t+h}) - \bar p_{t}(x_{t+h}) + \bar p_{t}(x_{t+h}) - \bar p_{t}(x_{t}),
\end{aligned}
\]
where,
\(
x_{t+h} = \Phi_{t, t+h}(x_t)\), \(\bar p_{t}(\mathrm x) =  \ell(\bar \Phi_{t+h, T}(\bar \Phi_{t, t+h}(\mathrm x))) = \bar p_{t+h}(\bar \Phi_{t,t+h}(\rm x))\) by the property (\ref{fami}).
Denote by $K \doteq \{x_t\colon t \in I\} \subset X$ the orbit of $x$, which is compact in $X$ as a continuous image of the compact interval $I$. Then, the expression under the sign of the desired limit is estimated by the sum of the terms $A^v_{t,h}$, $v \in \{\bar u, u\}$, 
\[
A^v_{t,h} \doteq \!\!\!\! \sup_{s \in I, {\mathrm x} \in K}\!\Big|\bar p_{s}(\Phi^v_{t,t+h}(\mathrm x)) - \bar p_{s}(\mathrm x) - h D \bar p_{s}(\mathrm x)\bigl[G({\mathrm{x}}) v_t\big] \Big|,
\]
 and \(D \bar p_t(x_{t+h})\Bigl[\big[G(x_{t+h})-G(x_t)\big]\bar u_t\Bigr]\). The letter tend to zero by continuity.
% we have 
% \[
% \begin{gathered}
% |A_{t,h}| \leq \sup_{s \in I, \, \mathrm x \in K} |\bar p_s(\bar \Phi_{t, t+h}(\mathrm x)) - \bar p_s(\mathrm x)|,\\ |B_{t,h}| \leq \sup_{s \in I, \, \mathrm x \in K} |\bar p_s(\Phi_{t, t+h}(\mathrm x)) - \bar p_s(\mathrm x)|. 
% \end{gathered}
% \]
Let us prove that $A^v_{t,h} = o_t(h)$.

2) \emph{Differentiability along the flow}. Leveraging the integral representation \eqref{Phi_integral_repr}, and adding and subtracting the term $ h F(\mathrm x, v_t)$, we can represent:
\begin{gather*}
\Phi^v_{t,t+h}({\rm x}) = {\rm{x}} +  h F(\mathrm x, v_t)+ \, \zeta_{t, \mathrm x}
(h),
\\
\zeta_{t, \mathrm x}(h) = \int_t^{t+h}\!\!\!\big(F(\Phi^v_{t, \tau}(\mathrm x), v_\tau) - F(\mathrm x, v_t)\big) \d \tau,  
\end{gather*}
and estimate ($f(\rm x)$ is killed in the last difference):
\begin{align*}
\|\zeta_{t, \mathrm x}(h)\|_X  & \leq  \int_{t}^{t+h} \Bigl\|F(\Phi^v_{t, \tau}(\mathrm x), v_\tau) - F(\mathrm x, v_\tau)\Bigr\|_X\d \tau \\
& + \int_t^{t+h} \!\!\bigl\|G(\mathrm x)[v_\tau - v_t]\big\|_X\d \tau.
\end{align*}
The first term is majorated by
\begin{gather*}
    \int_t^{t+h} \|\bar\Phi_{t, \tau} - {\id}_{X}\|_X \, (\lip_{K}(f)+ \lip_K(G) \|v_\tau\|_U) \d \tau\\ \leq M_K(\Phi)  h^2 (\lip_K(f)+ R\, \lip_K(G)) \d \tau,  
\end{gather*}
where we used the Lipschitz property: $ \|\bar\Phi_{t, \tau} - {\id}_{X}\|_X \leq M_K(\Phi) (\tau - t) \leq M_K(\Phi) h$ for $\tau \in [t, t+h]$. %and the Cauchy-Schwartz inequality $\|v\|_{L^1} \leq \sqrt T \, \|v\|_{L^2}$.  

The last is estimated by
\(
M_K(G) \int_t^{t+h} \|v_\tau - v_t\|_{U} \d \tau,  
\)
which is $o_t^v(h)$ according to the Lebesgue differentiation theorem \cite[II, Thm. 9]{diestel1977vector}. Since the final estimates are independent of $\mathrm{x}$, we conclude that $\sup_{\mathrm x \in {K}}\|\zeta_{t, \mathrm x}(h)\|_X \xrightarrow[]{h \downarrow 0}  0$.
Now, we leverage that fact that, for any fixed $s \in I$, \(\mathrm x \mapsto \bar p_s(\mathrm x) \doteq \ell(\Phi^v_{s, T}(\mathrm x))\) is continuously differentiable as a composition of $C^1$-functions. In view of the above analysis, we are then legal to write the $K$-uniform approximation{:}
\[
\sup_{K}\Big|\bar p_{s}(\Phi^v_{t,t+h}(\mathrm x)) - \bar p_{s}(\mathrm x) - h D \bar p_{s}(\mathrm x)\bigl[G({\mathrm{x}}) v_t\big] \Big| = o_{t,s}^v(h),
\]
in which
\(|o_{t,s}^v(h)| \leq \displaystyle\sup_{\rm x \in K}\|D\bar p_s(\mathrm x)\|_{\mathcal L}
\sup_{\mathrm x \in K}\|\zeta_{t, \mathrm x}(h)\|_X.\) 
The final step is to check that the first multiplier on the right is dominated by a constant, independent of $s$, so that $A^v_{t,h} = o_t(h)$ as desired.

3) \emph{Equi-differentiability{.}} The desired estimate follows from the definition of $\bar p$, yielding 
\[
\|D\bar p_s(\mathrm x)\|_{\mathcal L(X;\R)} \leq \sup_{\mathrm y \in K}\|D \ell(\mathrm y)\|_{\mathcal L(X;\R)} \|D\bar \Phi_{t, T}(\mathrm x)\|_{\mathcal L(X;X)},
\]
and the standard representation $D \bar \Phi_{t, \tau}(\mathrm x) = w^{t, \rm x}_\tau$ of the function $\tau \mapsto D \bar \Phi_{t, \tau}(\mathrm x)$ as a solution to the linear Cauchy problem in $\mathcal L(X;X)$ \cite[\S\S~~3.4-5]{cartan1983differential}:
\[
w_\tau^{t, \mathrm x} = {\id}_X + \int_t^\tau D F(\bar \Phi_{t, \tau}(\mathrm x), \bar u_\tau) \, w^{t, \mathrm x}_\tau \d \tau.
\]
By Gr\"{o}nwall's inequality, we get
\[
\|w_\tau^{t, \mathrm x}\|_{\mathcal L(X;X)} \leq \exp\big((M_{K}(Df) +  M_{K}(DG)\,  R) \, T\big). 
\]
Therefore,
\begin{align*}
&\sup_{\mathrm x \in K}\|D\bar p_s(\mathrm x)\|_{\mathcal L(X;\R)}\\[-0.1cm]
&\qquad\leq M_{K}(D\ell) \exp\big((M_K({Df}) +  M_K({DG})\, R) \, T\big).
\end{align*}
This observation completes the proof.

\bibliography{starmax_en_cleaned_sorted}

\begin{thebibliography}{16}
\providecommand{\natexlab}[1]{#1}
\providecommand{\url}[1]{\texttt{#1}}
\providecommand{\urlprefix}{URL }
\expandafter\ifx\csname urlstyle\endcsname\relax
  \providecommand{\doi}[1]{doi:\discretionary{}{}{}#1}\else
  \providecommand{\doi}{doi:\discretionary{}{}{}\begingroup
  \urlstyle{rm}\Url}\fi

\bibitem[{Amari(1977)}]{amariDynamicsPatternFormation1977}
Amari, S.i. (1977).
\newblock Dynamics of pattern formation in lateral-inhibition type neural
  fields.
\newblock \emph{Biological Cybernetics}, 27(2), 77--87.
\newblock \doi{10.1007/BF00337259}.

\bibitem[{Bensoussan(2007)}]{bensoussan2007representation}
Bensoussan, A. (2007).
\newblock \emph{Representation and Control of Infinite Dimensional Systems}.
\newblock Systems \& Control: Foundations \& Applications. Birkh{\"a}user
  Boston.

\bibitem[{Borz{\`\i}(2023)}]{borzi2023sequential}
Borz{\`\i}, A. (2023).
\newblock \emph{The Sequential Quadratic Hamiltonian Method: Solving Optimal
  Control Problems}.
\newblock Chapman \& Hall/CRC Numerical Analysis and Scientific Computing
  Series. CRC Press.

\bibitem[{Borz{\`\i} and Schulz(2011)}]{Borzi2011}
Borz{\`\i}, A. and Schulz, V. (2011).
\newblock \emph{{Computational Optimization of Systems Governed by Partial
  Differential Equations}}.
\newblock \doi{10.1137/1.9781611972054}.

\bibitem[{Cartan(1983)}]{cartan1983differential}
Cartan, H. (1983).
\newblock \emph{Differential Calculus}.
\newblock International studies in mathematics. Hermann.

\bibitem[{Diestel and Uhl(1977)}]{diestel1977vector}
Diestel, J. and Uhl, J. (1977).
\newblock \emph{Vector Measures}.
\newblock Mathematical surveys and monographs. American Mathematical Society,
  Providence, RI.

\bibitem[{Fattorini(1999)}]{Fattorini1999}
Fattorini, H.O. (1999).
\newblock \emph{Infinite Dimensional Optimization and Control Theory}.
\newblock Cambridge University Press.
\newblock \doi{10.1017/cbo9780511574795}.

\bibitem[{Kolokoltsov(2019)}]{kolokoltsov2019}
Kolokoltsov, V.N. (2019).
\newblock \emph{Differential equations on measures and functional spaces}.
\newblock Birkh\"{a}user, Cham.

\bibitem[{Krasovskii and Subbotin(2011)}]{krasovskii2011game}
Krasovskii, N. and Subbotin, A. (2011).
\newblock \emph{Game-theoretical control problems}.
\newblock Springer Series in Soviet Mathematics. Springer, New York.

\bibitem[{Krylov and Chernous’ko(1963)}]{Krylov1963}
Krylov, I. and Chernous’ko, F. (1963).
\newblock On a method of successive approximations for the solution of problems
  of optimal control.
\newblock \emph{USSR Computational Mathematics and Mathematical Physics}, 2(6),
  1371–1382.
\newblock \doi{10.1016/0041-5553(63)90353-7}.

\bibitem[{Li and Yong(1995)}]{Li1995}
Li, X. and Yong, J. (1995).
\newblock \emph{Optimal Control Theory for Infinite Dimensional Systems}.
\newblock Birkh\"{a}user Boston.
\newblock \doi{10.1007/978-1-4612-4260-4}.

\bibitem[{Pogodaev and Staritsyn(2024)}]{pogodaevExactFormulaeIncrement2024}
Pogodaev, N.I. and Staritsyn, M.V. (2024).
\newblock {Exact Formulae for the Increment of the Objective Functional and
  Necessary Optimality Conditions}, alternative to {{Pontryagin}}'s maximum
  principle.
\newblock \emph{Sbornik: Mathematics}, 215(6), 790--822.
\newblock \doi{10.4213/sm9967e}.

\bibitem[{Pontryagin et~al.(1962)Pontryagin, Boltyanskii, Gamkrelidze, and
  Mishchenko}]{Pontryagin1962}
Pontryagin, L.S., Boltyanskii, V.G., Gamkrelidze, R.V., and Mishchenko, E.F.
  (1962).
\newblock \emph{The Mathematical Theory of Optimal Processes}.
\newblock Interscience Publishers, a division of John Wiley \& Sons, New
  York-London.

\bibitem[{Srochko(1982)}]{srochko1982computational}
Srochko, V. (1982).
\newblock Computational methods of optimal control.
\newblock \emph{ISU, Irkutsk}.

\bibitem[{Staritsyn et~al.(2022)Staritsyn, Pogodaev, Chertovskih, and
  Pereira}]{SChPP-2022}
Staritsyn, M., Pogodaev, N., Chertovskih, R., and Pereira, F.L. (2022).
\newblock Feedback maximum principle for ensemble control of local continuity
  equations: An application to supervised machine learning.
\newblock \emph{IEEE Control Systems Letters}, 6, 1046--1051.
\newblock \doi{10.1109/LCSYS.2021.3089139}.

\bibitem[{Tr{\"o}ltzsch(2010)}]{Troeltzsch2010}
Tr{\"o}ltzsch, F. (2010).
\newblock \emph{Optimal Control of Partial Differential Equations: Theory,
  Methods and Applications}, volume 112 of \emph{Graduate Studies in
  Mathematics}.
\newblock American Mathematical Society.

\end{thebibliography}

\end{document}